\documentclass{article}
\usepackage[centertags]{amsmath}
\usepackage{amsfonts}
\usepackage{amssymb}
\usepackage{amsthm}
\usepackage{newlfont}
\usepackage{graphics}

\input xy
\xyoption{all}

\newcommand{\cP}{\mathcal{P}}
\newcommand{\bN}{\mathbb{N}}

\newcommand{\F}{\mathcal{F}}
\newcommand{\C}{\mathcal{C}}

\newcommand{\G}{\mathcal{G}}

\newtheorem{thm}{Theorem}
\newtheorem{prop}[thm]{Proposition}
\newtheorem{defi}[thm]{Definition}
\newtheorem{lemma}[thm]{Lemma}
\newcommand{\ten}{\otimes}
\newcommand{\comment}[1]{}

\begin{document}

\title{A forest formula for the antipode in incidence Hopf algebras}
\author{Hillary Einziger}

\maketitle

\begin{abstract}We present a new formula for the antipode of incidence Hopf algebras.  This formula is expressed as an alternating sum over forests.  First, we prove the formula for incidence Hopf algebras of families of lattices by exhibiting a map from chains of a lattice to forests.  Then, we extend the definition and present an analogous formula for the antipode of incidence Hopf algebras of families of posets.  We characterize those families for which our formula is cancellation-free.
\end{abstract}

\section{Introduction}

Many combinatorial Hopf algebras can be realized as incidence Hopf algebras of families of posets.  The antipode of incidence Hopf algebras, when expressed in terms of the canonical basis of indecomposable posets in the family, is generally highly non-trivial; i.e., the antipode is a sum with many terms. In \cite{mh89}, Haiman and Schmitt presented a closed formula for the antipode of the Fa\`{a} di Bruno Hopf algebra.  Their formula was expressed as a sum over trees (which we reinterpret as a sum over forests.)  They proved that this antipode formula is equivalent to Lagrange inversion.  In \cite{hf05}, Figueroa found a similar forest formula for the antipode of the incidence Hopf algebra of distributive lattices and presented applications in quantum field theory.  In this paper, we introduce a forest formula for the antipode of an arbitrary incidence Hopf algebra.  Both Figueroa's and Haiman and Schmitt's results are special cases of this new formula.

Haiman and Schmitt's formula for the antipode of the Fa\`{a} di Bruno Hopf algebra is cancellation-free, while Figueroa's formula for the antipode of the incidence Hopf algebra of distributive lattices is generally not cancellation-free.  We characterize those families of posets for which our forest formula for the antipode is cancellation-free.  The formula is cancellation-free for all indecomposable posets in a hereditary family if and only if every upper interval of every indecomposable interval in the family is indecomposable.  This condition is equivalent to the right-sided condition defined by Loday and Ronco in \cite{jl08}, which, in turn, is equivalent to the Lie algebra of primitive elements in fact being a pre-Lie algebra.

In Section 2, we recall the definitions of hereditary families of posets and incidence Hopf algebras, and we give the formula for the antipode as an alternating sum over chains.  In Section 3, we present several lemmas concerning the center and decomposition of posets.  In Section 4, we define forests of lattices.  In Section 5, we prove our forest formula for the antipode of incidence Hopf algebras of lattices.  In Section 6, we characterize those families for which the forest formula is cancellation-free.  In Section 7, we generalize our results from families of lattices to analogous results for incidence Hopf algebras of families of posets.

\section{Hereditary families}

An interval is a partially ordered set $P$ with unique
maximal and minimal elements, which we denote $\hat{1}_{P}$ and $\hat{0}_{P}$,
respectively. We eliminate the subscript when there is no chance of
ambiguity.  In this paper, we will assume all intervals are finite.

We slightly modify the definition of hereditary families of posets
from \cite{mh89}.

\begin{defi} \label{heredef}
A \emph{hereditary family} $(\cP, \sim, \cdot)$ is a family $\cP$
of finite intervals with a product operation $\cdot: \cP \times \cP
\rightarrow \cP$, where we write $PQ$ for $P \cdot Q$; and an
equivalence relation $\sim$ such that $\cP$ is closed under the
formation of subintervals and, for all $P, Q, R \in \cP$,
\begin{enumerate}
    \item If $P \sim Q$, then $PR \sim QR$.
    \item $(PQ)R \sim P(QR)$, and if $I$ is any single-element
    interval, then $IP \sim PI \sim P$.  Also, $PQ \sim QP$.
    \item If $P \sim Q$, then there is a bijection $x \mapsto x'$ from
    $P$ to $Q$ such that $[\hat{0}_{P},x] \sim [\hat{0}_{Q},x']$ and $[x, \hat{1}_{P}] \sim [x', \hat{1}_{Q}]$ for
    all $x \in P$.
    \item There is a poset isomorphism $\psi$ from the Cartesian
    product $P \times Q$ to $PQ$ such that $[\psi(\hat{0}_{P}, \hat{0}_{Q}), \psi(x,y)] \sim [\hat{0}_{P},x][\hat{0}_{Q},y]$
    and $[\psi(x,y), \psi(\hat{1}_{P}, \hat{1}_{Q})] \sim [x, \hat{1}_{P}][y, \hat{1}_{Q}]$ for all $(x,y) \in P \times Q$.
\end{enumerate}
We use $\prod$ to denote the iterated $\cdot$ product.
\end{defi}

Let $k$ be a commutative ring with $1$, and let $\cP = (\cP, \sim,
\cdot )$ be a hereditary family. Conditions (1) and (2) of Definition \ref{heredef} imply that the
quotient $\cP /\! \sim$ is a commutative monoid with product induced
by the product in $\cP$. The \textit{incidence Hopf algebra} of
$\cP$, denoted $H(\cP)$, is the monoid algebra of $\cP /\! \sim$ over
$k$, with coalgebra structure given by
$$\delta(P) = \sum_{x \in P} [\hat{0},x] \ten [\hat{x},1],$$
and
$$
    \epsilon(P) =
    \begin{cases}
        1   &\text{if $|P| = 1$}\\
        0   &\text{otherwise.}
    \end{cases}
$$

This coproduct is clearly coassociative.  Condition (3) ensures
that $\delta$ is well-defined on $\cP /\! \sim$.  Condition (4)
guarantees that $\delta(PQ) \sim \delta(P) \delta(Q)$, and so
$H(\cP)$ is a bialgebra.

We do not distinguish notationally between elements of $\cP$ and
$\cP /\! \sim$.

\begin{defi} Let $\cP$ be a hereditary family of posets.  We say
that a poset is \emph{decomposable} in $\cP$ if it is the
non-trivial product of non-singleton posets in $\cP$, and we say that it is
\emph{indecomposable} otherwise.  If $P \in \cP$, then we write $D(P)$ for
the set of all $x \in P \setminus \{\hat{0}, \hat{1}\}$ such that $[\hat{0}, x]$ is
decomposable in $\cP$, and we write $I(P)$ for the set of all $x \in P
\setminus \{\hat{0}, \hat{1}\}$ such that $[\hat{0},x]$ is indecomposable in $\cP$.
\end{defi}

Let $\cP_{0}$ be the set of all indecomposable posets in $\cP$.
Then $\cP_0 /\! \sim$ is the free commutative monoid on $\cP /\! \sim$.
So, as an algebra, $H(\cP)$ is isomorphic to the polynomial algebra
$k[\cP_{0} /\! \sim]$.

A \textit{chain} $C$ in an interval $P$ is a set $ \hat{0} = c_{0} <
c_{1} < \dots < c_{n} = \hat{1}$ of elements of $P$, and $\ell(C) = n$ is
the length of $C$. The \textit{length} of an interval $P$ is the
length of its longest chain.  Let $H(\cP)_{n}$ be the submodule of $H(\cP)$ spanned by intervals
of length less than or equal to $n$.  These define a bialgebra
filtration $H(\cP)_{0} \subseteq H(\cP)_{1} \subseteq H(\cP)_{2}
\subseteq \dots$.  Condition (2) of Definition \ref{heredef} ensures that $H(\cP)_{0} = k \cdot
1$, and so $H(\cP)$ is a connected bialgebra, and thus it is a Hopf
algebra.

\begin{defi} The \emph{convolution algebra} of $H(\cP)$ is $Hom_{k}\big(H(\cP), H(\cP)\big)$ with operation
convolution defined by
$$(f \ast g)(P) = \sum_{x \in P} f\big([\hat{0},x]\big)g\big([x,\hat{1}]\big).$$
\end{defi}

Since $H(\cP)$ is a Hopf algebra, the identity map $id: H(\cP)
\rightarrow H(\cP)$ in the convolution algebra has a two-sided
convolution inverse $\chi$, known as the \textit{antipode}.

For $P \in \cP$, let $\C(P)$ be the set of all chains of $P$. For
each $P \in \cP$, define $\Omega: \C(P) \rightarrow \cP$ by

$$\Omega (C) = \prod_{i=1}^{\ell(C)} [c_{i-1},
c_{i}].$$

\begin{prop}\label{chaindef} \cite{mh89}
The antipode of $H(\cP)$ is given by
$$\chi(P) = \sum_{C \in \C(P)} {(-1)}^{\ell(C)}
\Omega (C).$$
\end{prop}

\section{The center of a poset}

In this section, we consider $\cP$ to be a hereditary family of
posets.

\begin{defi} Let $P \in \cP$.  The \emph{center} of $P$,
denoted $Z(P)$, is the set of all $a \in P$ such that there is some
$a' \in P$ such that $[\hat{0},a][\hat{0},a'] \sim [\hat{0},\hat{1}]$.  The \emph{prime center} of
$P$, denoted $Z'(P)$, is the set of all minimal non-zero elements of
$Z(P)$.
\end{defi}

In \cite{gb67}, Birkhoff described several properties of the centers
of posets.  In his work, two posets were considered to be equivalent
if they were isomorphic, and the product considered was Cartesian
product.  We show that these properties hold in the more general case of hereditary families of posets.

For $P \in \cP$, let $I_{\cong}(P)$ be the set of all $x
\in P$ such that $[\hat{0},x]$ is indecomposable as a Cartesian product.
Let $D_{\cong}(P)$ be the set of all $x \in P$ such that $[\hat{0},x]$ is
decomposable as a Cartesian product, and let $Z_{\cong}(P)$ be the
set of all $a \in P$ such that there is some $a' \in P$ such that
$[\hat{0},a][\hat{0},a'] \cong [\hat{0},\hat{1}]$, where the product is Cartesian product.

It is clear from Condition (4) of the hereditary family definition that $D(P) \subseteq D_{\cong}(P)$ and $I_{\cong}(P) \subset I(P)$ for any $P \in \cP$.  It then follows that $Z(P) \subseteq Z_{\cong}(P)$.  However, it is not necessarily true that $Z'(P) \subseteq Z'_{\cong}(P)$, since the minimal elements of $Z(P)$ may not be minimal in the larger set $Z_{\cong}(P)$.

\begin{lemma}\label{jmcenter} If $a \in Z(P)$, where $P \in \cP$, then $a \vee z$ and $a \wedge z$ exist
for any $z \in P$.  Additionally, if $a \in Z(P)$, then $P \sim
[\hat{0},a][a,\hat{1}]$ by the map $z \rightarrow (z \wedge a, z \vee a)$, where
$z \in P$.
\end{lemma}

\begin{proof}
Birkhoff shows the first assertion for any $a \in Z_{\cong}(P)$, and so it is true for any $a \in Z(P)$.

The second assertion is a modification of
a lemma from Birkhoff.  Since $a \in Z(P)$, we know that $P \sim XY$ where $X \sim [\hat{0},a]$
and $Y$ is some poset.  So, in the $\psi$ map described in
Condition (4) of Definition \ref{heredef}, we have $a =
\psi (\hat{1}_{X}, \hat{0}_{Y})$.  Let $z \in P$.  Then $z = \psi (x,y)$ for some $x \in X$ and $y \in Y$.  Then, by the given map,

$$z = \psi (x,y) \rightarrow \big(\psi (x,y) \wedge \psi(\hat{1}_{X},\hat{0}_{Y}), \psi(x,y) \vee \psi(\hat{1}_{X},\hat{0}_{Y})\big) = \big(\psi(x,\hat{0}_{Y}),
\psi(\hat{1}_{X},y)\big).$$

The elements of the form $\psi(\hat{1}_{X},y)$ in the factorization $P \sim
XY = [\hat{0},a]Y$ are exactly the elements of $[a,\hat{1}]$ in $P$.  So then we must have
$Y \sim [a,\hat{1}]$, and so the given map sends $P$ to $[\hat{0},a][a,\hat{1}]$, as
needed.
\end{proof}

\begin{lemma}\label{compcenter}If $P \in \cP$ and $b \in Z(P)$, then $b$ has a
unique complement $b' \in Z(P)$, and $[\hat{0},b][\hat{0},b'] \sim P$.
\end{lemma}

\begin{proof}
Birkhoff shows the existence of the unique complement in $Z_{\cong}(P)$, and the existence of the unique complement in $Z(P)$ follows from the same reasoning. The
previous lemma shows that $P \sim [\hat{0},b][b,\hat{1}]$ by the map $z
\rightarrow (z \wedge b, z \vee b)$.  Then $b' \rightarrow (\hat{0},\hat{1})$,
and so we must have $[\hat{0},b'] \sim [b,\hat{1}]$, and thus $P \sim
[\hat{0},b][\hat{0},b']$.
\end{proof}

Birkhoff uses the analogues of the previous two lemmas to prove the following three lemmas when the equivalence relation is isomorphism and the product is Cartesian product.  The general hereditary family case follows from the same reasoning.

\begin{lemma} $Z(P)$ is a Boolean lattice and a sublattice of $P$.
\end{lemma}

\begin{lemma} If $a \in Z(P)$, then $a \in Z'(P)$ if and only if $a \in I(P)$.
\end{lemma}

\begin{lemma}\label{prodcenter} $P \sim \prod_{a \in Z'(P)} [\hat{0},a]$.
\end{lemma}

\begin{defi}
An element $a$ of a lattice $P$ is said to be \emph{distributive}
if the identities
$$a \wedge (x \vee y) = (a \wedge x) \vee (a \wedge y)$$
$$x \wedge (a \vee y) = (x \wedge a) \vee (x \wedge y)$$
and their duals hold for all $x, y \in P$.  An element $a$ is
\emph{complemented} in $P$ if there exists an element $a' \in P$
such that $a \vee a' = \hat{1}$ and $a \wedge a' = \hat{0}$.
\end{defi}

Birkhoff proves the following lemma.

\begin{lemma}\label{distcenter}
If $a \in P$, then $a \in Z_{\cong}(P)$ if and
only if $a$ is both distributive and complemented in $P$.
\end{lemma}

\section{Forests of lattices}

For the next several sections, we consider $\cP$ to be a family of lattices.  The more general poset case is considered in Section 7.

\begin{defi}\label{forestdef}
A \emph{forest} of a lattice $P \in \cP$ is a set $F \subseteq
I(P)$, with $\bigvee F \neq \hat{1}$ and $\hat{0} \notin F$, such that:
\begin{enumerate}
    \item if $a_{1}$, $a_{2} \in F$, then either $a_{1} \leq a_{2}$,
    $a_{2} \leq a_{1}$, or $a_{1} \wedge a_{2} = \hat{0}$. (This condition is referred to as ``non-overlapping.")
    \item if $\{b_{i}\}_i$ is an antichain in $F$, then
    $\prod_{i}[\hat{0}, b_{i}] \sim [\hat{0}, \bigvee_{i} b_{i}]$.
\end{enumerate}
\end{defi}

\noindent \textbf{Example 1.} In \cite{hf05}, Figueroa defined forests of distributive lattices.  Figueroa's definition relied on the fundamental theorem of distributive lattices: If $L$ is a finite distributive lattice, then $L$ is isomorphic to the poset of order ideals $J_{P}$ of some finite poset $P$.  He then defined a forest $F$ of $L$ as a collection of connected order ideals of $P$ where $\emptyset \notin F$ and $\bigcup F \neq P$, such that if $I_{1}, I_{2} \in F$, then either $I_{1} \cap I_{2} = \emptyset$, or $I_{1} \subseteq I_{2}$, or $I_{2} \subseteq I_{1}$.  In the map sending order ideals of $P$ to $J_{P}$, connected order ideals are sent to $I(J_{P})$, and so this non-overlapping condition on order ideals is equivalent to Condition 1 of Definition \ref{forestdef}.  The second condition of Definition \ref{forestdef} holds for any antichain in a distributive lattice.

\noindent \textbf{Example 2.} In the lattice shown in Figure 1, the forests of the interval are the empty
forest; the single-element forests $\{a\}$ and $\{b\}$; and the
two-element forest $\{a,b\}$.
Note that $c$ cannot be in any forest since the interval $[\hat{0},c]$ is
decomposable.

\begin{figure}[h]
$$
\hbox{ \begin{picture}(120,150)
    \put(60,30){\circle*{3}}
    \put(30,60){\circle*{3}} \put(90,60){\circle*{3}}
    \put(60,90){\circle*{3}}
    \put(60,120){\circle*{3}}
    \put(60,30){\line(-1,1){30}} \put(60,30){\line(1,1){30}}
    \put(30,60){\line(1,1){30}} \put(90,60){\line(-1,1){30}}
    \put(60,90){\line(0,1){30}}
    \put(25,60){$\scriptstyle a$} \put(95,60){$\scriptstyle b$}
    \put(55,90){$\scriptstyle c$}
    \end{picture}
}$$
\caption{Example 2}
\end{figure}
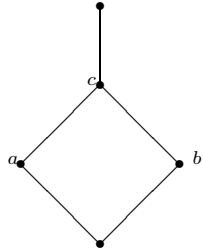

\noindent \textbf{Example 3.} In the lattice shown in Figure 2, the forests of the interval are the empty
forest; the single-element forests $\{a\}$, $\{b\}$, $\{c\}$, and
$\{d\}$; and the three two-element forests $\{a,d\}$, $\{b,d\}$, and
$\{c,d\}$.  The set $\{a, b\}$, for example, is not a forest,
because even though it satisfies the first condition of the forest definition, it does not satisfy the second condition, since $a \vee b = d$
and $[\hat{0}, a][\hat{0},b] \nsim [\hat{0}, d]$.

\begin{figure}[h]
$$
\hbox{ \begin{picture}(120,150)
    \put(60,30){\circle*{3}}
    \put(30,60){\circle*{3}} \put(90,60){\circle*{3}}
    \put(60,60){\circle*{3}}
    \put(60,90){\circle*{3}}
    \put(60,120){\circle*{3}}
    \put(60,30){\line(-1,1){30}} \put(60,30){\line(1,1){30}}
    \put(60,30){\line(0,1){30}}
    \put(30,60){\line(1,1){30}} \put(90,60){\line(-1,1){30}}
    \put(60,60){\line(0,1){30}}
    \put(60,90){\line(0,1){30}}
    \put(25,60){$\scriptstyle a$} \put(95,60){$\scriptstyle c$}
    \put(55,60){$\scriptstyle b$}
    \put(55,90){$\scriptstyle d$}
    \end{picture}
}$$
\caption{Example 3}
\end{figure}
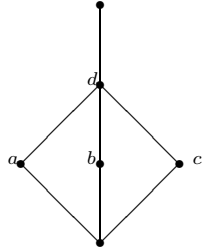

\section{Forest formula for antipode in incidence Hopf algebras of hereditary families of
lattices}\label{forestsection}
If $F$ is a forest of a lattice $P$ and $b \in F \cup \{\hat{1}\}$, then
we say $a$ is a \textit{predecessor} of $b$ in $F$ if $a \in F$,
$a < b$, and there is no $a' \in F$ such that $a < a' <
b$. If $b$ has no predecessors in $F$, then we consider $\hat{0}$ to be
its predecessor.

\begin{defi}\label{phidef}If $P \in \cP$, then let $\F(P)$ be the set
of all forests of $P$.  If $F \in \F(P)$, then let
$$\Theta(F) = \prod_{b \in F \cup Z'(P)}
[\tilde{b}, b]$$
where $\tilde{b}$ is the join of all the predecessors of $b$ in
$F$.  Note that
$$\Theta(F) = \prod_{b \in F \cup \{\hat{1}\}} [\tilde{b}, b]$$
is an equivalent definition.
\end{defi}

\begin{prop}\label{forestsurj}
Let $P \in \cP$ be a finite lattice.  There is a surjection $\phi:
\C(P) \rightarrow \F(P)$, with $C \mapsto F_C$, such that
$\Omega(C) \sim \Theta(F_{C})$ for all $C \in \C(P)$.
\end{prop}

\begin{proof}

If $C: \hat{0} = c_{0} < c_{1} < \dots < c_{\ell(C)} = \hat{1}$ is a chain of $P$, then let

$$F_{C} := \bigcup_{i=1}^{\ell(C)-1} Z'([\hat{0}, c_{i}]).$$

We want to show that $F_{C}$ satisfies Definition \ref{forestdef}.  We know that
$[\hat{0}, a]$ is indecomposable for any $a \in F_{C}$, since $a$ is in the prime center of some $[\hat{0}, c_i]$.  Also, we know that $\bigvee F_{C} = c_{\ell(C) - 1} < \hat{1}_{P}$, and so $F_{C}$ satisfies the preliminary conditions of Definition \ref{forestdef}.

    Now, we need to show that $F_{C}$ satisfies the first condition of Definition \ref{forestdef}. Let $a,b \in F_{C}$.  We want to show that that either
    $a \leq b$, $b \leq a$, or $a \wedge b = \hat{0}$.  Let $i$ and $j$ be
    the smallest indices such that $a$ and $b$ are in the prime
    centers of $[\hat{0}, c_{i}]$ and $[\hat{0}, c_{j}]$, respectively.  Without
    loss of generality, assume $i \leq j$.  We know that $a \leq
    c_{i} \leq c_{j}$.

    If $a \leq b$, then we are done.  If not, then we must have $b <
    c_{j}$.  Since $b \in Z'([\hat{0}, c_{j}])$, there is a unique $b' <
    c_{j}$ such that $b \wedge b' = \hat{0}$ and $[\hat{0},c_{j}] \sim [\hat{0},b][\hat{0},b']$.
     Since $[\hat{0},a]$ is indecomposable and $a < c_{j}$, we know that
     either $a \leq b$ or $a \leq b'$.  By assumption, the first
     possibility is false, and so we must have $a \leq b'$, and thus
     $a \wedge b = \hat{0}$.

    Next, we want to show that, if $\{b_{i}\}_i$ is an antichain in
    $F_{C}$, then each $b_{i}$ is in the center of $[\hat{0}, \bigvee
    b_{i}]$.  We induct on the size of the antichain.

    First, suppose we have an antichain $\{b_{1}, b_{2}\}$ in
    $F_{C}$.  As before, let $i$ and $j$ be the smallest indices
    such that $b_{1}$ and $b_{2}$ are in the prime centers of $[\hat{0},
    c_{i}]$ and $[\hat{0}, c_{j}]$, respectively.  Assume $i \leq j$.

    Since $b_{2}$ is in the prime center of $[\hat{0}, c_{j}]$, we know from Lemma \ref{compcenter}
    that there is some $b' < c_{j}$ such that $[\hat{0}, b_{2}][\hat{0}, b']
    \sim [\hat{0}, c_{j}]$.  We know that $b_{1} \leq c_{i} \leq c_{j}$,
    and, since we showed that $F_C$ satisfies the first forest condition, we know that $b_{1} \wedge b_{2} =
    \hat{0}$, and so we must have $b_{1} \leq b'$.  Then by Condition (4)
    of Definition \ref{heredef}, we must have $[\hat{0}, b_{1}][\hat{0},
    b_{2}] \sim [\hat{0}, b_{1} \vee b_{2}]$.

    Now assume that $\prod [\hat{0}, b_{i}] \sim [\hat{0}, \bigvee b_{i}]$ for all
     antichains in $F_{C}$ up to size $n-1$.  Let $\{b_{1}, \dots , b_{n}\}$ be an antichain in
     $F_{C}$.  For each $1 \leq i \leq n$, let $k_{i}$ be the
     smallest index such that $b_{i} \in Z'([\hat{0}, c_{k_{i}}])$.
     Without loss of generality, assume $k_{n} \geq k_{i}$ for all
     $i$.  Since $b_{n}$ is in the center of $[\hat{0}, c_{k_{n}}]$, we know from Lemma \ref{distcenter} that $b_n$ is
     distributive in $[\hat{0}, c_{k_n}]$, and so the distributive
     property shows
     $$b_{n} \wedge (\bigvee_{j=1}^{n-1} b_{j}) =
     \bigvee_{j=1}^{n-1} (b_{n} \wedge b_{j}) = \hat{0},$$
     since, by the previous part, $b_{n} \wedge b_{j} = \hat{0}$ for any
     $1 \leq j \leq n-1$.

     Since $b_{n}$ is in the center of $[\hat{0}, c_{k_{n}}]$, Lemma \ref{compcenter} shows that it must have a unique complement $b'$ in $[\hat{0}, c_{k_{n}}]$ such that $[\hat{0}, b_{n}][\hat{0}, b'] \sim [\hat{0},
     c_{k_{n}}]$.  We know that $\bigvee_{j=1}^{n-1} b_{j}
     \leq c_{k_{n}}$ and we just showed that $b_{n} \wedge
     (\bigvee_{j=1}^{n-1} b_{j}) = \hat{0}$, and so we must have
     $\bigvee_{j=1}^{n-1} b_{j} \leq b'$.  So then, by
     induction and Condition (4) of the hereditary family
     definition, we get
     $$\prod_{j=1}^{n} [\hat{0}, b_{j}] = [\hat{0}, b_{n}] \prod_{j=1}^{n-1}
     [\hat{0}, b_{j}] \sim [\hat{0}, b_{n}][\hat{0}, \bigvee_{j=1}^{n-1} b_{j}] \sim
     [\hat{0}, b_{n} \vee \bigvee_{j=1}^{n-1} b_{j}] = [\hat{0},
     \bigvee_{j=1}^{n} b_{j}],$$
     as needed.

Next, we want to show that this map is surjective.  Let $F$ be a
forest of $P$. Let $S_{1}$ be the collection of all maximal elements
of the forest, $S_{2}$ the collection of all maximal elements of $F
\setminus S_{1}$, and so forth.  Define $C_F$ as the chain $\hat{1}
> \bigvee S_{1} > \bigvee S_{2}
> \dots > \hat{0}$.  We need to show that $\bigvee S_{i} > \bigvee S_{i+1}$ for all $i$.

Clearly, $\bigvee S_{i} \geq \bigvee S_{i+1}$.  Also, $Z'([\hat{0}, \bigvee S_{i}]) =
S_{i}$ and $Z'([\hat{0}, \bigvee S_{i+1}]) = S_{i+1}$.  Since $S_{i}$ and $S_{i+1}$
are disjoint, we know that $\bigvee S_{i} \neq \bigvee S_{i+1}$, and so $\bigvee S_{i} >
\bigvee S_{i+1}$. So each forest is, in fact, associated to at least one
chain.

Last, we want to show that $\Omega(C) \sim \Theta(F_{C})$ for any
chain $C$ of $P$. Let $\ell(C) = n$.  For each $1 \leq i \leq n$,
let $F_{i} = Z'([\hat{0}, c_{i}])$. Then $\bigcup_{i} F_{i} = F_{C}
\cup Z'(P)$.  For each $i$, we know from Lemma \ref{prodcenter} that

$$[\hat{0}, c_{i}] \quad \sim \prod_{a \in Z'([\hat{0}, c_{i}])} [\hat{0},a].$$

Then, from Lemma \ref{prodcenter} and Condition (4) of the hereditary family definition, we know that

$$[\hat{0}, c_{i-1}] \quad \sim \prod_{a \in Z'([\hat{0}, c_{i}])} [\hat{0}, a \wedge c_{i-1}].$$

We then get

$$[c_{i-1}, c_{i}] \quad \sim \prod_{a \in Z'([\hat{0}, c_{i}])} [a \wedge c_{i-1}, a].$$

For $a \in Z'([\hat{0}, c_{i}]$, we find

$$ a \wedge c_{i-1} = a \wedge \bigg(\bigvee_{b \in Z'([\hat{0}, c_{i-1}])} \!\!\!b\bigg)\quad = \bigvee_{b \in Z'([\hat{0}, c_{i-1}])}\!\!\!\! a \wedge b\quad = \bigvee_{\substack{b \in Z'([\hat{0}, c_{i-1}])\\ b \leq a}} \!\!\!\!b \quad = \quad\tilde{a}.$$

So then $[c_{i-1}, c_{i}] \sim \prod_{a \in F_{i}} [\tilde{a}, a]$, and so we get $\Omega(C) \sim \Theta(F_{C})$.

\end{proof}

The first main result of this paper is a Zimmerman-type formula for
the antipode of $H(\cP)$ in terms of forests.  To derive that
formula, we require one more proposition.

\begin{prop}\label{forestsum}If $F$ is a forest of $P$, then $\sum_{C \in {\phi}^{-1}(F)} (-1)^{\ell(C)} = (-1)^{d(F)}$.
\end{prop}

\begin{proof}

     We follow the similar argument of \cite{mh89}.  Let a
     filtration $G$ of the forest $F$ be a chain $\emptyset = I_{0}
     \subset I_{1} \subset \dots \subset I_{k} = F$ of lower order ideals of
     $F$ such that, for all $1 \leq j \leq k$, the set $I_{j} \setminus I_{j-1}$
     is an antichain.  The length of the filtration is $\ell(G) = k$.  We
     claim that, for each forest $F$ of a poset $P$, there is a
     bijection between ${\phi}^{-1}(F)$ and the set $\G(F)$ of all filtrations
     of $F$, such that, if a chain $C$ is mapped to the filtration $G$,
     then $\ell(C) = \ell(G) + 1$.

     First, suppose $C \in {\phi}^{-1}(F)$ is $\hat{0} = c_{0} < c_{1} < \dots < c_{n}
     = \hat{1}$.  Define the filtration $G$ by setting $I_{0} = \emptyset$
     and $I_{k} = Z'([\hat{0}, c_{k}]) \cup I_{k-1}$ for $1 \leq k \leq
     n-1$.  Each $I_{k}$ must be a lower order ideal of $F$, and
     the $I_{k}$ must be strictly increasing. Since $I_{k} \setminus I_{k-1}$ is a subset of $Z'([\hat{0},
     c_{k}])$, it must be an antichain.  Thus, $G$ is a filtration.

     Conversely, given a filtration $\emptyset = I_{0} \subset I_{1}
     \subset \dots \subset I_{n} = F$ of $F$, define the chain $C$ by letting $c_{0} = \hat{0}$ and
     $c_{k} = \bigvee I_{k}$ for $1 \leq k \leq n$, and let $c_{n+1} = \hat{1}$. Then,
     by our definition of a forest, $Z'([\hat{0}, c_{k}])$ must be the set of
     maximal elements of $I_{k}$.  Each element of $I_{k} \setminus
     I_{k-1}$ must be either greater than some maximal element of
     $I_{k-1}$ or not comparable to any element of $I_{k-1}$, and so
     we have $c_{k} > c_{k-1}$.

     These constructions are clearly inverse to one another, and so they
     form a bijection.

    Lemma 4 of \cite{mh89} states that, if $Q$ is a finite
    poset and $\G(Q)$ is the set of all filtrations of $Q$, then
    $$\sum_{G \in \G(Q)} (-1)^{\ell(G)} = (-1)^{|Q|}.$$

     By using this lemma and the given bijection, we get
     $$ \sum_{C \in {\phi}^{-1}(F)} (-1)^{\ell(C)} = \sum_{G \in \G(F)}
     (-1)^{\ell(G) + 1} = (-1)^{d(F)}.$$

\end{proof}

We now come to our first main result.

\begin{thm}\label{forestantth}
If $\cP$ is a hereditary family of lattices, then the antipode of $H(\cP)$ is given by
\begin{equation}\label{forestant}
\chi(P) = \sum_{F \in \F(P)}(-1)^{d(F)}\Theta(F)
\end{equation}
for all $P \in \cP$.
\end{thm}

\begin{proof}

We know from Proposition \ref{chaindef} that
$$\chi(P) = \sum_{C \in \C(P)} (-1)^{\ell(C)} \Omega(C).$$

Using Proposition \ref{forestsurj}, we get
$$\chi(P) = \sum_{F \in \F(P)} \Theta(F) \sum_{C \in {\phi}^{-1}(F)}
(-1)^{\ell(C)}.$$

Then Proposition \ref{forestsum} gives us
$$\chi(P) = \sum_{F \in \F(P)} (-1)^{d(F)} \Theta(F).$$

\end{proof}

\section{Conditions for non-cancellation in computation of antipode}

Formula \eqref{forestant} is very similar to the formula of Zimmerman,
explored in \cite{hf05}, for the antipode of the Hopf algebra of
Feynman graphs.  Zimmerman's antipode formula has the useful
property of being cancellation-free.

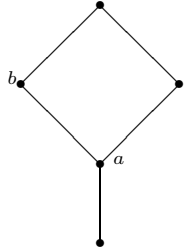
\begin{figure}[h]
    $$
    \hbox{ \begin{picture}(120,150)
    \put(60,30){\circle*{3}} \put(30,90){\circle*{3}}
    \put(60,60){\circle*{3}} \put(90,90){\circle*{3}}
    \put(60,120){\circle*{3}}
    \put(60,30){\line(0,1){30}} \put(60,60){\line(-1,1){30}}
    \put(60,60){\line(1,1){30}} \put(30,90){\line(1,1){30}}
    \put(90,90){\line(-1,1){30}}
    \put(65,60){$\scriptstyle a$} \put(25,90){$\scriptstyle b$}
    \end{picture}
    }$$
    \caption{Lattice without cancellation-free forest computation}\label{Fi:cancellations}
\end{figure}

For the forest computation of $\chi(P)$ given by \eqref{forestant} to be
cancellation-free, it must be the case that $(-1)^{d(F)}$ and $(-1)^{d(F')}$ have the same sign whenever $F$ and $F'$ are forests of $P$ such that $\Theta(F) \sim \Theta(F')$.  In general, Formula \eqref{forestant} is not cancellation-free. A simple
example of an indecomposable lattice for which the forest
computation of $\chi$ is not cancellation-free is shown in Figure
\ref{Fi:cancellations}. If $P$ is the lattice shown, then the
forests $F = \{a\}$ and $F' = \{a,b\}$ will cancel each other in
the computation of $\chi(P)$.

It is also clear that the forest computation of $\chi(P)$ will have
cancellations if $P$ is decomposable: if $a \in Z'(P)$, then
$\Theta(\{a\}) = [\hat{0},a][a,\hat{1}] \sim P = \Theta(\emptyset)$.  Note that, since $\chi$ is multiplicative, it is determined by its value on indecomposables.  We now characterize those indecomposable lattices in hereditary families for which the forest computation of the antipode is cancellation-free.  We then characterize those hereditary families for which the forest computation is cancellation-free for all indecomposable lattices.

\begin{defi} An indecomposable lattice $P \in \cP$ is called \emph{upper-indecomposable} if, for every $x<\hat{1} \in P$, the interval $[x,\hat{1}]$ is indecomposable. An indecomposable lattice $P \in \cP$ is called \emph{super-upper-indecomposable} (s.u.i.) if every indecomposable interval of $P$ is upper-indecomposable.
\end{defi}

\begin{prop}A lattice $P \in \cP$ is s.u.i if and only if every indecomposable lower interval of $P$ is upper-indecomposable in $\cP$.
\end{prop}

\begin{proof}

Clearly, if $P$ is s.u.i., then every indecomposable lower interval of $P$ is upper-indecomposable.

We prove the converse by contradiction.  We want to
    show that if every indecomposable lower interval of $P$ is upper-indecomposable, and there are $a, x, y \in P$ such that $\hat{0} \leq a
    < x < y \leq \hat{1}$ and $[x,y]$ is decomposable, then $[a,y]$
    is decomposable.  So suppose $a$, $x$, and $y$ are as given, and
    $[x,y]$ is decomposable.  Then, since every indecomposable lower interval of $P$ is upper-indecomposable, the interval
    $[\hat{0},y]$ must be decomposable.  Let $Z'([\hat{0},y]) =
    \{y_{i}\}_{i=1}^{n}$.

    Let $x_{i} = x \wedge y_{i}$ for all $i$.  Then Condition (4) of the hereditary family definition implies
 $[\hat{0},x] \sim \prod_{i=1}^{n} [\hat{0}, x_{i}]$ and $[x,y]
    \sim \prod_{i=1}^{n} [x_{i}, y_{i}]$.  Since $[\hat{0}, y_{i}]$ is
    indecomposable for each $i$, and every indecomposable lower interval of $P$ is upper-indecomposable, we know that $[x_{i}, y_{i}]$ is
    indecomposable for each $i$.  Since $[x,y]$ is decomposable, we
    must have $x_{i} \lneq y_{i}$ for at least two values of $i$.

    Let $a_{i} = a \wedge y_{i}$ for all $i$.  Since $a < x$, we must also have $a_{i} \leq x_{i}$ for all $i$.
    Since $x_{i} \lneq y_{i}$ for at least two values of $i$, we must have $a_{i} \lneq y_{i}$ for at least two values
    of $i$, and so, since $[a,y] \sim \prod_{i=1}^{n} [a_{i},
    y_{i}]$, we conclude that $[a,y]$ is decomposable, as needed.

\end{proof}

\begin{prop}
Let $P \in \cP$ be an indecomposable lattice.  The forest computation of
$\chi(P)$ given by Theorem \ref{forestantth} is cancellation-free if and only if $P$ is s.u.i.
\end{prop}

\begin{proof}
    Suppose an indecomposable lattice $P$ is s.u.i.  Let $F$ be a forest
    of $P$.  Then $\Theta(F) = \prod_{b \in F \cup
    \{\hat{1}\}} [\tilde{b}, b]$.  Since $[\hat{0}, b]$ is indecomposable for each
    $b$ and $P$ is s.u.i., we know that each $[\hat{0},b]$ is upper-indecomposable, and so each $[\tilde{b}, b]$ will be
    indecomposable.  We also know that $\tilde{b} \neq b$.  So then $d(F)$ is the number of indecomposable intervals in the unique factorization of $\Theta(F)$. Thus, the
    computation of $\chi(P) = \sum_{ F \in \F(P)} (-1)^{d(F)}
    \Theta(F)$ will be cancellation-free.

    Conversely, suppose $P$ is indecomposable but not s.u.i. Then
    there are some $x, y \in P$ such that $\hat{0} < y < x$ and $[\hat{0}, x]$ is
    indecomposable, but $[y, x]$ is decomposable.  Say $[y, x] \sim [y,
    x_{1}] \cdots [y, x_{n}]$ is the factorization of $[y,x]$ into indecomposable intervals.

    Assume $x < \hat{1}$.  (The case where $x = \hat{1}$ is analogous.)

    Case 1: Suppose $[\hat{0}, x_{i}]$ is indecomposable for some $i$.
    Let $F = Z'([\hat{0},y]) \cup \{x\}$, and let $F' = Z'([\hat{0},y]) \cup
    \{x_{i}, x\}$.  Then
    $$\Theta (F) = [\hat{0},y][y,x][x,\hat{1}],$$
    and
    $$\Theta (F') = [\hat{0},y][y,x_{1}][x_{1},x][x,\hat{1}] \sim
    [\hat{0},y][y,x][x,\hat{1}].$$

    So $\Theta (F) \sim \Theta (F')$, but $d(F') = d(F) + 1$, so
    these two forests cancel each other in the forest computation of
    $\chi(P)$.

    Case 2: Suppose $[\hat{0}, x_{i}]$ is decomposable for all $i$.  Let
    $$[\hat{0}, x_{1}] \sim [\hat{0}, x_{1,1}][\hat{0}, x_{1,2}] \cdots [\hat{0},
    x_{1,k}]$$
    be the unique factorization of $[\hat{0}, x_1]$ into indecomposables.
    We know that $[y, x_{1}]$ is indecomposable, and so, without
    loss of generality, we can say $y = (y_{1,1}, x_{1,1}, x_{1,2}, \dots , x_{1,k})$
    in this factorization of $[\hat{0}, x_{1}]$.  We can see that
    $[y_{1,1}, x_{1,1}] \sim [y, x_{1}]$.  We can also see that
    $Z'([\hat{0},y]) = Z'([\hat{0}, y_{1,1}]) \cup \{x_{1,2}, \dots x_{1,k}\}$.

    Let $F = Z'([\hat{0},y]) \cup \{x\}$ and let $F' = Z'([\hat{0},y]) \cup
    \{x_{1,1}, x\}$.  These both satisfy Definition \ref{forestdef}, and we have
    $$\Theta (F) = [\hat{0},y][y,x][x,\hat{1}],$$
    and
    $$\Theta (F') = [\hat{0},y][y_{1,1}, x_{1,1}][x_{1}, x][x,\hat{1}] \sim
    [\hat{0},y][y, x_{1}][x_{1},x][x,\hat{1}] \sim [\hat{0},y][y,x][x,\hat{1}].$$

    So again, $\Theta (F) \sim \Theta (F')$, but $d(F') = d(F) +
    1$, and hence these forests cancel each other.

\end{proof}

\begin{defi} A hereditary family $\cP$ is called \emph{upper-indecomposable} if every indecomposable $P \in \cP$ is upper-indecomposable.
\end{defi}

Note that, since a hereditary family must be closed under the taking of intervals, the hereditary family $\cP$ is upper-indecomposable if and only if every indecomposable $P \in \cP$ is super-upper-indecomposable.

These propositions bring us to our next main result.

\begin{thm}
Let $\cP$ be a hereditary family.  Then the forest computation of
$\chi(P)$ given by Theorem \ref{forestantth} will be cancellation-free for all indecomposable $P \in
\cP$ if and only if $\cP$ is upper-indecomposable.
\end{thm}

In \cite{jl08}, Loday and Ronco defined a \emph{cofree-coassociative combinatorial Hopf algebra} as a cofree bialgebra $H$ together with an isomorphism between $H$ and the tensor coalgebra over the primitive elements of $H$.  Furthermore, such $H$ satisfies the \emph{right sided condition} if $\delta(Q(H)) \subseteq H \ten Q(H)$, where $Q(H)$ denotes the subspace of irreducibles in $H$.  The upper-indecomposable hereditary families $\cP$ are exactly those hereditary families for which $H(\cP)$ satisfies the right-sided condition.  Theorem 5.3 of \cite{jl08} states that the right-sided cofree-coassociative combinatorial Hopf algebras are exactly those cofree-coassociative combinatorial Hopf algebras in which the primitive elements form a pre-Lie algebra, rather than merely a Lie algebra.

\noindent \textbf{Example 4.} The partition lattice $\Pi_{n}$ of the set $\{1, \dots , n\}$ is the poset of all partitions of $\{1, \dots , n\}$, ordered by refinement: if $x, y \in \Pi_{n}$, then $x \leq y$ if every block of $x$ is contained in a block of $y$.  The Fa\`{a} di Bruno Hopf algebra is the Hopf algebra $H(\cP)$,
where $\cP$ is the set of all finite products of finite partition lattices, with the equivalence relation $\sim$ given by isomorphism.  In \cite{mh89}, Haiman and Schmitt define a surjection from the set of chains of $\Pi_{n}$ to the set of leaf-labelled trees with $n$ leaves and no vertices of degree $1$.  According to their definition, if $C$ is a chain in the partition lattice $\Pi_{n}$, then the tree $T(C)$ associated with $C$ is the poset of all subsets of $\{1, \dots , n\}$ which
appear as blocks of partitions in $C$, ordered by inclusion.

If $x$ is a partition in $\Pi_{n}$ with non-singleton blocks $B_{1}, \dots , B_{k}$, then $Z'([\hat{0},x])$ is the set of all partitions $\{a_{1}, \dots , a_{k}\}$ of $\Pi_{n}$, where $a_{i}$ is the partition with block $B_{i}$ and all other elements as singleton blocks.  For a chain $C$, the forest $\phi(C)$ given by Proposition \ref{forestsurj} is the set of all partitions in $\Pi_{n}$ with one non-singleton block, such that the non-singleton block is a block of a partition in $C$.  For each chain $C$ of $\Pi_{n}$, then, there is a clear bijection between the forest $\phi(C)$ and the tree $T(C)$, since the non-singleton blocks represented by the internal vertices of $T(C)$ are exactly the non-singleton blocks in the elements of the forest $\phi(C)$.  The refinement order on $\Pi_{n}$ inherited by $\phi(C)$ is the same as the ordering of the blocks in $T(C)$ by inclusion.

The indecomposable $P \in \cP$ are the members of the equivalence classes of the partition lattices $\Pi_{n}$ for all $n$.  Each upper interval $[\rho, \hat{1}]$ in a partition lattice is equivalent to $\Pi_{|\rho|}$, and so $\cP$ is an upper-indecomposable family, and thus the forest computation of $\chi(P)$ is cancellation-free for all $\Pi_{n}$.  Haiman and Schmitt proved that the Lagrange inversion formula is equivalent to this antipode formula.

\section{Forest formula for the antipode for hereditary families of
posets}

The antipode formula in Theorem \ref{forestantth} for hereditary families of lattices can be extended to a
formula for the antipode of the Hopf algebra of any hereditary
family of posets.
Clearly, the second condition of Definition \ref{forestdef} cannot be
applied to general posets, since general posets lack a join
operation. A poset $P$ with non-overlapping indecomposable lower
intervals $[\hat{0},a]$ and $[\hat{0},b]$ might have several elements $c > a,b$ such that $[\hat{0},a][\hat{0},b] \sim [\hat{0},c]$.  Note that the converse, however, is
not true.  Any decomposable lower interval of a poset has a unique
factorization as a product of indecomposable lower intervals.

We introduce a new definition of forest.

\begin{defi}\label{posetfor} A \emph{forest} of a poset $P \in \cP$ is an ordered pair $(F, J_{F})$ where $F \subseteq I(P) \setminus \{\hat{0}, \hat{1}\}$ and $J_{F} : 2^{F} \rightarrow P \setminus \{\hat{1}\}$ such that:
\begin{enumerate}
    \item If $a_{1}, a_{2} \in F$, then either $a_{1} \leq a_{2}$; $a_{2} \leq a_{1}$; or $[\hat{0}, a_{1}] \cap [\hat{0}, a_{2}] = \{\hat{0}\}$.
     \item If $G \subseteq F$ and $\{a_{i}\}$ is the set of maximal elements of $G$, then $[\hat{0}, J_{F}(G)] \sim \prod [\hat{0}, a_{i}]$.
     \item If $G' \subseteq G \subseteq F$, then $J_{F}(G') \leq J_{F}(G).$
\end{enumerate}
\end{defi}

Note that a single set $F$ can have several different $J_{F}$ functions which each satisfy this definition.  Each of these $(F, J_{F})$ pairs is regarded as a unique forest.  If $P$ is a lattice, however, then Conditions (2) and (3) of Definition \ref{posetfor} guarantee that $J_{F}$ must be the join operation, and so Definition \ref{forestdef} can be seen as a special case of Definition \ref{posetfor}.

Note that Conditions (1) and (2) of Definition \ref{posetfor} are essentially the same as the conditions of Definition \ref{forestdef}.  The join operation of a lattice always satisfies Condition (3) of Definition \ref{posetfor}. It is possible, however, for a non-lattice finite interval $P$ to
    have a subset $F$ and a function $J: 2^{F} \rightarrow P$
    such that $F$ and $J$ satisfy the first two
    conditions of Definition \ref{posetfor}, but not the third, and so the third condition is not superfluous.

\subsection{Motivating example: the $N$-colored Fa\`{a} di Bruno Hopf algebra}

We generalize the Fa\`{a} di Bruno Hopf algebra described in Example 4.  Following the example of \cite{mh89}, let $N \in \bN$.  An
$N$-colored set is a finite set $X$ with a map $\theta = \theta_{X}:X
\rightarrow \{1, \dots, N\}$, with $\theta(x)$ called the color of
$x$.  Let $X_{r} = \{x \in X | \theta (x) = r\}$, and if $X$ is an
$N$-colored set, let $|X|$ be the vector $(|X_{1}|, \dots ,
|X_{N}|)$.  An $N$-colored partition of an $N$-colored set $X$ is a
partition $\pi$ of $X$ such that each block of $\pi$ is assigned a
color and, if $\{x\}$ is a singleton block of $\pi$, then
$\theta_{\pi}(\{x\}) = \theta_{X}(x)$.

The poset $\Pi_{\textbf{n}}$ of $N$-colored partitions of an
$N$-colored set $X$ with $|X| = \textbf{n} = (n_{1}, \dots , n_{N})$ is formed by letting
$\pi \leq \rho$ if $\pi \leq \rho$ in the refinement order of the partition lattice and, if $B$ is a
block of both $\pi$ and $\rho$, then $\theta_{\pi}(B) =
\theta_{\rho}(B)$.  This poset has a $\hat{0}$ but it does not have a $\hat{1}$; the partition with a single block can occur with any of
$N$ colors.  Let $\Pi_{\textbf{n}}^{r}$ be the poset of $N$-colored
partitions with all maximal elements except the one colored $r$
deleted.

Let $\mathfrak{F}^{N}$ be the family of all $\Pi_{\textbf{n}}^{r}$ of $N$-colored
partitions posets of $N$-colored sets.  If $\pi \leq \rho$ in some $N$-colored partition poset, the let $\rho | \pi$ be the $N$-colored partition of the set $\pi$ induced by the colors of the unions of the blocks of $\pi$ in $\rho$. Define the relation $\sim$ as
color-isomorphism: $[\pi, \rho] \sim [\pi', \rho']$ if $[\pi, \rho]
\cong [\pi', \rho']$ and there is a bijection $B \mapsto B'$ from
the non-singleton blocks of $\rho | \pi$ to the non-singleton blocks
of $\rho' | \pi'$ such that
\begin{enumerate}
    \item If $B$ is the union of $i$ blocks of $\pi$, then $B'$ is
    the union of $i$ blocks of $\pi'$;
    \item $\theta_{\rho | \pi} (B) = \theta_{\rho' | \pi'} (B')$;
    \item For each non-singleton block $B$ of $\rho | \pi$, there is
    a bijection $\varphi$ from the blocks of $\pi$ in $B$ to the
    blocks of $\pi'$ in $B'$ such that, if $\beta$ is a block of
    $\pi$ in $B$, then $\theta_{\pi}(\beta) =
    \theta_{\pi'}(\varphi(\beta))$.
\end{enumerate}

For example, if we let subscripts denote the color of an element or
block, then
$$[1_{1}/2_{1}/3_{2}/4_{2}/5_{1}, \ (13)_{2}/(45)_{1}/2_{1}]$$
$$\sim [1_{1}/2_{1}/3_{2}/(45)_{2} , \ (145)_{2}/(23)_{1}]$$
$$\sim [1_{1}/2_{2}/3_{1}/4_{2}, \ (12)_{2}/(34)_{1}].$$

Let the product on the equivalence classes be defined so that
$$[\pi, \rho] \sim \prod_{B \in \rho | \pi}
\Pi_{|B|}^{\theta_{\rho | \pi}(B)}.$$

For example,
$$[1_{1}/2_{2}/3_{1}, \ (12)_{1}/ 3_{1}][1_{2}/(23)_{1}, \
(123)_{2}]$$
$$\sim [1_{1}/2_{2}/3_{1}/4_{2}/(56)_{1}, \
(12)_{1}/3_{1}/(456)_{2}]$$
$$\sim [1_{1}/2_{2}/3_{2}/4_{1}, \ (12)_{1}/(34)_{2}].$$

It is straightforward to verify that the family of all finite products of the elements of $\mathfrak{F}^{N}$ is a
hereditary family.  The incidence algebra of $\mathfrak{F}^{N}$ is
known as the \emph{$N$-colored Fa\`{a} di Bruno Hopf algebra}.

As in the Fa\`{a} di Bruno Hopf algebra, if $P = \Pi_{\textbf{n}}^{r} \in \mathfrak{F}^{N}$, then
$I(P)$ will be the set of all $\pi \in P$ such that $\pi$
has exactly one non-singleton block. Although the posets in $\mathfrak{F}^{N}$ are not lattices, there is a unique choice of $J_{F}$ for each possible forest set $F$.  If $P = \Pi_{\textbf{n}}^{r}$, then, as in the Fa\`{a} di Bruno Hopf algebra,
    a forest $(F, J_{F})$ must be a set of partitions in $P$ such that each
    partition has exactly one non-singleton block and, if $\pi, \rho \in
    F$ such that $B_{\pi}, \ B_{\rho}$ are their respective
    non-singleton blocks, then either $B_{\pi} \subseteq B_{\rho}$,
    $B_{\rho} \subseteq B_{\pi}$, or $B_{\pi} \cap B_{\rho} =
    \emptyset$.  The only possible map $J_{F}$ which satisfies Condition $(2)$ of Definition \ref{posetfor} is that in which $J_{F}(G)$ is the $N$-colored partition whose non-singleton blocks are exactly the non-singleton blocks of the maximal elements of $G$, with the same colors.

\subsection{Forest antipode formula for hereditary families of posets}

 Using arguments similar to the proofs in Section \ref{forestsection}, we find a forest formula for the antipode of incidence Hopf algebras of hereditary families of posets.

\begin{defi} Let $P \in \cP$, let $(F, J_{F})$ be a forest of $P$, and let $a \in F$.  Then let $\tilde{a} = J_{F}(\{p(a)\})$, where $\{p(a)\}$ is the set of all predecessors of $a$ in $F$.  As in Definition \ref{phidef}, we let $\Theta(F, J_F ) = \prod_{a \in F \cup Z'(P)} [\tilde{a}, a]$
\end{defi}

\begin{thm} \label{posetant}
Let $\cP$ be a hereditary family of posets.  A formula for the antipode in $H(\cP)$ is
$$\chi(P) = \sum_{(F, J_F) \in \F(P)}(-1)^{d(F)}\Theta(F, J_F )$$
for all $P \in \cP$.
\end{thm}

\begin{proof}

First, as in Proposition \ref{forestsurj} we find a surjection $\phi: \C(P) \rightarrow \F(P)$.  As in Proposition \ref{forestsurj}, if $C$ is a chain of $P$, then let
$$F_{C} = \phi(C) := \bigcup_{i=1}^{\ell(C)-1} Z'([\hat{0}, c_{i}]).$$

Next, we define $J_{F_{C}}$ by induction.  First, $J_{F_{C}}(\emptyset) = \hat{0}$, and, if $a \in F_{C}$, then $J_{F_{C}}(\{a\}) = a$.

Now, let $\{b_{i}\}_{i=1}^{n}$ be an antichain in $F_{C}$.  For each $1 \leq i \leq n$, let $c_{k_{i}}$ be the minimal chain element such that $b_{i} \in Z'([\hat{0}, c_{k_{i}}])$.  Assume $k_{n} \geq k_{i}$ for all $i$.  Assume by induction that we have defined $J_{F_{C}} (\{b_{i}\}_{i=1}^{n-1})$, and that $J_{F_{C}} (\{b_{i}\}_{i=1}^{n-1}) \leq c_{k_{n-1}} \leq c_{k_{n}}$.  Since $b_{n} \in Z'([\hat{0}, c_{k_{n}}])$, we know by Lemma \ref{jmcenter} that $b_{n} \vee J_{F_{C}} (\{b_{i}\}_{i=1}^{n-1})$ is defined in $[\hat{0}, c_{k_{n}}]$, and so let $J_{F_{C}} (\{b_{i}\}_{i=1}^{n}) = b_{n} \vee J_{F_{C}} (\{b_{i}\}_{i=1}^{n-1})$ in $[\hat{0}, c_{k_{n}}]$.  Since join is commutative and associative, $J_{F_{C}} (\{b_{i}\}_{i=1}^{n})$ is well-defined.

We can regard the elements of $F_{C}$, together with the image of $J_{F_{C}}$, as a subposet of $P$.  In fact, if $G$ is a subset of $F_{C}$, then $J_{F_{C}}(G)$ is a minimal upper bound of $G$.  We can thus regard $F_{C}$ together with the image of $J_{F_{C}}$ as a ``sublattice" of $P$, in which joins in the sublattice correspond to $J_{F_{C}}$ in $P$.  The proofs of Propositions \ref{forestsurj} and \ref{forestsum} and Theorem \ref{forestantth} can thus be easily modified from the lattice case to the general poset case, completing the proof.
\end{proof}

The conditions for non-cancellation in a hereditary family of lattices similarly generalize to families of posets.

\begin{thm}The antipode calculation in Theorem \ref{posetant} is cancellation-free for all indecomposable $P$ in a hereditary family of posets $\cP$ if and only if $\cP$ is upper-indecomposable.
\end{thm}

\noindent \textbf{Example 5}
In \cite{fc07}, Chapoton and Livernet define a hereditary family of posets from an set-operad $\cP$.  If $\cP$ is a set-operad, then $\Pi_{\cP}$ is the species Comm$\circ \cP$, where Comm is the species that maps a finite set $I$ to the singleton $\{I\}$.  For each finite set $I$, they introduce a partial order on $\Pi_{\cP}(I)$.  Proposition 3.4 of \cite{fc07} proves that, for each set-operad $\cP$, the family of all $\Pi_{\cP}(I)$ is a hereditary family.  Proposition 3.3 shows that this family satisfies our upper-indecomposable condition.  Therefore, the computation of $\chi$ given by Theorem \ref{posetant} is cancellation-free for the incidence Hopf algebra of this family.


\begin{thebibliography}{99}

\bibitem{gb67}
    Garett Birkhoff,
    \emph{Lattice Theory},
    Amer. Math. Soc. Colloq. Publ. XXV, Providence, R.I., (1967).
\bibitem{fc07}
    F. Chapoton and M. Livernet,
    \emph{Relating two Hopf algebras built from an operad}
    arXiv:0707.3725v1.
\bibitem{hf05}
    H\'{e}ctor Figueroa,
    \emph{Combinatorial Hopf algebras in quantum field thoery I},
    Rev.Math.Phys., 17, (2005), 881.
\bibitem{mh89}
    Mark Haiman and William Schmitt,
    \emph{Incidence Algebra Antipodes and Lagrange Inversion in One
    and Several Variables},
    Combinatorial Theory, 50, (1989), 172.
\bibitem{jl08}
    Jean-Louis Loday and Mar\`{i}a Ronco,
    \emph{Combinatorial Hopf algebras},
    Clay Math. Proc., 10, (2008).
\end{thebibliography}
\end{document}